\documentclass[a4j,12pt]{article}
\usepackage{amsmath,amsthm,amssymb,bm}

\numberwithin{equation}{section}

\newtheorem{thm}{Theorem}[section]
\newtheorem{lem}[thm]{Lemma}
\newtheorem{prop}[thm]{Proposition}
\newtheorem{cor}[thm]{Corollary}
\newtheorem{definition}[thm]{Definition}
\newtheorem{remark}[thm]{Remark}
\newtheorem{example}[thm]{Example}

\newenvironment{rem}{\begin{remark}\rm}{\end{remark}}
\newenvironment{ex}{\begin{example}\rm}{\end{example}}

\title{Around Shannon's Interpretation for Entropy-preserving Stochastic Averages}

\author{Marie Choda}

\date{}

\begin{document}
\maketitle
\centerline{Osaka Kyoiku University, Osaka, Japan} 
\centerline{marie@cc.osaka-kyoiku.ac.jp}
\begin{abstract} 
Let $\Phi$ be a positive unital  Tr-preserving map  on $M_n(\mathbb {C})$.  
We give various characterizations for $\Phi$ to preserve the  von Neumann entropy of a state $\rho$ on $M_n(\mathbb {C})$.  
Among others, it is given by that  $\Phi$ behaves to $\rho$ as a *-automorphism.  
This  is also  equivalent to that the entropy of the stochastic matrix arising  from $\{\rho, \Phi \}$ is  zero. 

\end{abstract} 

keywords: {Entropy, state, positive linear map, bistochastic matrix}

\smallskip

\section{Introduction} 
In this note, from the view point of the von Neumann entropy,  
we study the pair $\{\rho,  \Phi \} $ of a state $\rho$ and a 
positive unital \text{Tr}-preserving map $\Phi$ on the algebra 
$M_n(\mathbb {C})$  of $n \times n$ complex matrices, 
where  $\text{Tr}$ is  the standard trace of $M_n(\mathbb {C})$ so that $\text{Tr} (e) = 1$ 
for every minimal projection $e$ in $M_n(\mathbb {C})$. 

For a state $\rho$ of $M_n(\mathbb {C})$, the von Neumann entropy $S(\rho)$ is given as  the von Neumann entropy $S(D_\rho)$ 
of the density operator $D_\rho$ of $\rho$.  
If $\Phi$ is a positive unital \text{Tr}-preserving map on $M_n(\mathbb {C})$, then 
$\Phi(D_\rho)$ plays the  density operator of a new state of $M_n(\mathbb {C})$ that is $\rho \circ \Phi^*$. 
It holds in general that  $S(\rho) \le S(\rho \circ \Phi^*)$ 
(i.e. $S(D_\rho) \le S(\Phi(D_\rho))$) (see for example \cite{NS, OP, Pe}). 

We give in Theorem 3.3 below some characterizations for the equality of this relation. 
For an example, 
$S(D_\rho) = S(\Phi(D_\rho))$ if and only if  $\Phi$ behaves for $\rho$ as a *-automorphism of $M_n(\mathbb {C})$, i.e. 
$\rho \circ \Phi^* = \rho \circ \alpha$ for some *-automorphism $\alpha$ of $M_n(\mathbb {C})$, in other words 
$\Phi(D_\rho) = u D_\rho u^*$ for some unitary $u$ in $M_n(\mathbb {C})$. 

This our formulation is an extended version of Shannon's interpretation in  \cite [p.395 4.]{Sch} 
(a detailed description of which we denote in Section 3.1) 
for  entropy-preserving stochastic averages 
\vskip 0.3cm

As a main tool, we pick up  a bistochastic matrix $b_\rho(\Phi)$ arising from a given pair $\{\rho, \Phi\}$ of a state $\rho$ and a positive unital 
Tr-preserving map $\Phi$. 

Based on the notion of weighted entropy for a bistochastic matrix introduced in \cite{ZSKS}, 
we define the entropy $H^\lambda(b_\rho(\Phi))$ of this bistochastic matrix $b_\rho(\Phi)$ with 
respect to the probability vector $\lambda$  by the eigenvalues of $D_\rho$. 

Another entropic object $S_\rho(\Phi)$ (which we give below) is arising as 
an average of a family of von Neumann entropy coming from  $D_\rho$ via $\Phi$.  
\smallskip 

We show in Theorem 3.10 the following relations:   
\begin{align*}
S_\rho(\Phi) \leq H^\lambda(b_\rho(\Phi)) &\le S(\rho \circ \Phi^*)    \le  S(\rho) + S_\rho(\Phi) \\
 &\hskip 1.1cm \vee \Vert \\
 &\hskip 1.1cm   S(\rho)
\end{align*}
together with  characterizations for the each  equality, and in  Theorem 3.14 we show that 
$S(D_\rho) = S(\Phi(D_\rho))$ if and only if $H^\lambda(b_\rho(\Phi)) = 0$. 
\vskip 0.3cm

\section{Preliminaries}
Here we summarize notations, terminologies and basic facts. 
\vskip 0.3cm

\noindent{\bf 2.1.}  
The {\it entropy function} $\eta$ is defined on  $[0,1]$ by 
$$\eta(t) = -t \log t \quad  (0 < t \leq 1) \quad \text{and} \quad \eta(0) = 0.$$
The $\eta$ is {\it strictly operator-concave},  i.e.  
for a $k$-tuple of real numbers $\{t_i \}_{i=1}^k$ such that $t_i > 0, \sum_{i=1}^k t_i = 1$ 
and bounded self-adjoint operators $\{x_i\}_{i=1}^k$ with spectra in [0,1], 
it  holds in general that 
$$\sum_{i=1}^k t_i \eta(x_i) \le \eta( \sum_{i=1}^k t_i x_i)$$
and the equality implies that $x_i = x_j$ for all $i, j$. 
(see for example \cite [B]{NS}, \cite{OP, Pe}). 
\vskip 0.3cm

\noindent{\bf 2.2.} 
Let $\lambda = ( \lambda_1, \cdots, \lambda_n )$ be a probability vector in $\mathbb {R}^n$.  
The  {\it Shannon entropy} $H(\lambda)$ for  $\lambda$  is given as 
$$H(\lambda) = \eta(\lambda_1) + \cdots + \eta(\lambda_n).$$ 
It holds always that $H(\lambda) \le \log n$ 
and $H(\lambda) \leq H(\lambda  b)$ (\cite{Sch}, cf.\cite{PP}) 
for a {\it bistochastic} matrix  $b = [b_{ij}]$ 
(i.e. $b_{ij} \geq 0, \ \sum_i b_{ij} = \sum_j b_{ij} = 1$ for all $i,j = 1, \cdots, n$). 
A bistochastic matrix  is also called a {\it doubly stochastic}  matrix (see for example  \cite{Pe}).  
\vskip 0.3cm

\noindent{\bf 2.3.} 
Every positive linear functional $\phi$ on $M_n(\mathbb {C})$ is of the form 
$\phi(x) = \text{Tr}(D_\phi x )$, $(x  \in M_n(\mathbb {C}))$ for a unique positive element $D_\phi$ in $M_n(\mathbb {C})$, 
which is called the {\it density operator} or {\it density matrix}. 
The density matrix  of a state $\rho$ is  characterized as a positive operator $D_\rho$ such that $\text{Tr}(D_\rho) = 1$. 

Using the eigenvalue list $\{\lambda_1,  \cdots, \lambda_n \}$  of $D_\phi$, 
the {\it von Neumann entropy} $S(\phi)$ for a positive linear functional 
$\phi$ and the {\it von Neumann entropy} $S(D_\phi)$  for $D_\phi$ are defined  by 
$$S(\phi) = S(D_\phi) = \sum_{i=1}^n \eta(\lambda_i). $$
If $\rho$ is a state of $M_n(\mathbb {C})$, then the eigenvalue vector $\lambda = (\lambda_1,  \cdots, \lambda_n)$  of $D_\rho$ is 
a probability vector and so $S(\rho) = S(D_\rho) = H(\lambda)$. 
\smallskip

\section{Pairs $\{\rho, \Phi\}$ of  states  and   positive maps} 
Let $\rho$ be a  state of $M_n(\mathbb {C})$, and let $D_\rho$ be the density matrix of $\rho$. 
Let $\Phi : M_n(\mathbb {C}) \rightarrow  M_n(\mathbb {C})$ be a positive unital \text{Tr}-preserving map. 
Then $\Phi(D_\rho)$ is a positive operator in $M_n(\mathbb {C})$ and $\rm{Tr}(\Phi(D_\rho)) =  1$.  

In order to see the  state whose density matrix is $\Phi(D_\rho)$, we need the system of the Hilbert-Schmidt inner product 
of  $M_n(\mathbb {C})$. 
The inner product is given by 
$<x, y> \ = \ \text{Tr} (y^*x)$ for $x,y \in M_n(\mathbb {C})$ and 
the adjoint map $\Phi^* : M_n(\mathbb {C}) \rightarrow  M_n(\mathbb {C})$ of $\Phi$ 
is given by 
$\text{Tr} (y \Phi^*(x)) = \text{Tr} (\Phi(y) x  )$ for $x, y \in M_n(\mathbb {C})$. 

Since $\Phi$ is positive and Tr-preserving,  it implies that $\Phi^*$ is positive and 
unital so that $\rho \circ \Phi^*$ is  a state, whose density matrix is $\Phi(D_\rho)$:  
$$\rho  \circ \Phi^* (x) = \text{Tr} (D_\rho \Phi^*(x) ) = \text{Tr} (\Phi(D_\rho) x), \quad (x \in  M_n(\mathbb {C})).$$
\vskip 0.3cm

Let 
\begin{equation*}
\lambda = (\lambda_1, \cdots, \lambda_n) \quad \text{and} \quad \mu = (\mu_1, \cdots, \mu_n),
\end{equation*} 
be the  probability vectors of the eigenvalues of $D_\rho$ and $\Phi(D_\rho)$
respectively.  
Here we arrange them always in a decereasing order: 
\begin{equation*}
\lambda_1 \geq \lambda_2 \geq  \cdots \geq \lambda_n  \quad \text{and} \quad \mu_1 \geq \mu_2 \geq  \cdots \geq \mu_n. 
\end{equation*}

Let $\{e_1, \cdots, e_n\}$  (resp. $\{p_1, \cdots, p_n\}$) be the 
mutually orthogonal minimal projections inducing the  following decomposition of  $D_\rho$  (resp. $\Phi(D_\rho)$):  
\begin{equation}
D_\rho = \sum_{i=1}^n \lambda_i e_i \quad \text{and} \quad \Phi(D_\rho) = \sum_{j=1}^n \mu_j p_j.
\end{equation}

We denote by $u_{(\rho, \Phi)}$  the unitary such that  
\begin{equation}
u_{(\rho, \Phi)} e_i {u_{(\rho, \Phi)}}^*  = p_i, \quad \text{for each} \quad i = 1, \cdots, n. 
\end{equation}

Also we denote by $A$ (resp. $B$) the maximal abelian subalgebra of $M_n(\mathbb {C})$ 
generated by $\{e_1, \cdots, e_n\}$ (resp. $\{p_1, \cdots, p_n\}$). 

\smallskip

\subsection{A generalization of Shannon's interpretation}
The motivation of this note is to give a generalized version of Shannon's interpretation for  
entropy-preserving stochastic averages of probability vectors 
in the framework of von Neumann entropy for states on $M_n(\mathbb {C})$. 

Shannon denotes in \cite [p.395 4.]{Sch}  as the followings: 
If we perform any "averaging " operation on the $p = \{p_i\}_{i=1,\cdots,n}$ of the form 
$$p_i' = \sum_j a_{ij} p_j$$
(where  $a_{ij} \geq 0, \ \sum_i a_{ij} = \sum_j a_{ij} = 1$), the entropy $H$ increases (except in the special case 
where this transformation amounts to no more than a permutation of the $p_i$ with $H$ of course remaining the same). 

A detailed explanation for this was presented in \cite{PP} together with some characterizations for  
pairs $\{p, [a_{ij}] \}$ of entropy-preserving stochastic averages.    

By replacing a probability vector  (resp. a bistachastic matrix) to a state $\rho$ of $M_n(\mathbb {C})$
(resp. a  unital positive Tr-preserving map $\Phi$ on $M_n(\mathbb {C})$), 
we show here that 
the action of $\Phi$ on $\rho$ preserves   the von Neumann entropy 
if and only if  $\Phi$ behaves just as an automorphism for the state $\rho$. 
\vskip 0.3cm

\smallskip
\subsubsection{Bistochastic matrix $b_\rho(\Phi)$ for the pair $\{\rho, \Phi\}$} 
A matrix $b = [b_{ij}] \in M_n(\mathbb {C})$ is called a bistochastic matrix if 
$$b_{ij}  \ge 0 \quad \text{and} \quad  \sum_i b_\rho(\Phi)_{ij} = \sum_j b_\rho(\Phi)_{ij} = 1, \ \text{for all} \ i,j$$ 

For a  state $\rho$ of $M_n(\mathbb {C})$  and 
a unital positive {\rm Tr}-preserving map $\Phi$  on $M_n(\mathbb {C})$, 
we denote by $b_\rho(\Phi)$ the matrix  given by 
\begin{equation}
b_\rho(\Phi)_{ij} = \text{Tr}(\Phi(e_i)p_j), \quad  (1 \leq i \leq n, 1 \leq j \leq n) . 
\end{equation}
\smallskip

Assume that $\lambda_i > \lambda_{i+1}$ and $\mu_i > \mu_{i+1}$ for all $i = 1, \cdots, n-1$. 
Then  the value $b_\rho(\Phi)_{ij}$ is uniquely 
determined for all $i,j$ because the spectral projections $\{e_i\}_i$ and $\{p_j\}_j$ are uniquely determined. 
In the other case, the value $b_\rho(\Phi)_{ij}$ is not always uniquely determined. 
For an example, if it happened that  
$D = \sum_i \lambda_i f_i$ and  $\Phi(D) =  \sum_j \mu_j q_i$ 
for some projections $\{f_i\}_{i=1}^n$ and $\{q_j\}_{j=1}^n$ different from $\{e_i\}_{i=1}^n$ and $\{p_j\}_{j=1}^n$,  
then the matrix $b_\rho(\Phi)^{e, p}$ with the  $i,j$-coefficient  $\text{Tr}(\Phi(e_i)p_j)$ may be different from 
the matrix $b_\rho(\Phi)^{f, q}$ with the  $i,j$-coefficient $\text{Tr}(\Phi(f_i)q_j)$.  
But the difference is covered by the permutation matrices  $[\sigma]$ and $[\pi]$ 
via the  permutations $\sigma$ and $\pi$  such that  $f_i = e_{\sigma(i)}$ and $q_i = p_{\sigma(i)}$ for all $i$ and $j$:  
$b_\rho(\Phi)^{f, q} = [\sigma] b_\rho(\Phi)^{e, p} [\pi]$. 

We show later that our disccusions do not depend on these kind of matrix representations and so we denote simply by $b_\rho(\Phi)$. 
\smallskip

\begin{lem} 

{\rm (1)} The $b_\rho(\Phi)$ is a bistochastic matrix. 

{\rm (2)} The $b_\rho(\Phi)$ transposes the vector $\lambda$ to the vector $\mu$, i.e.  $\lambda b_\rho(\Phi) = \mu.$
\end{lem}
\vskip 0.3cm

\noindent{\it Proof.} 
(1) 
Since $\Phi$ is  positive, it holds that  
$b_\rho(\Phi)_{ij} = \text{Tr}(p_j \Phi(e_i)p_j) \ge 0$ for all $i$ and $j$. 
The condition that $\text{Tr} \circ \Phi  = \text{Tr}$  implies that  for all $i$
$$\sum_j b_\rho(\Phi)_{ij} = \sum_j \text{Tr}(\Phi(e_i)  p_j) = 
\text{Tr}(\Phi(e_i) \sum_j p_j) = \text{Tr}(\Phi(e_i) ) = 1,$$ 
and condition that $\Phi (1_{M_n(\mathbb {C})}) = 1_{M_n(\mathbb {C})}$ implies that for all $j$ 
$$\sum_i b_\rho(\Phi)_{ij} = \sum_i \text{Tr}(\Phi(e_i)  p_j)  = \text{Tr}(\Phi(\sum_i e_i)  p_j) = \text{Tr}( p_j) = 1,$$ 
so that $b_\rho(\Phi)$ is a bistochastic matrix.  
\vskip 0.3cm

(2)  For all $j$, we have that 
$$\mu_j = \text{Tr}(\Phi(D_\rho) p_j) = \text{Tr}( \sum_{i=1}^n \lambda_i  \Phi(e_i) p_j) = \sum_{i=1}^n \lambda_i b_\rho(\Phi)_{ij}.$$
This means that $\mu = \lambda b_\rho(\Phi)$. 
\qed
\vskip 0.3cm

For each $j$, we set 
\begin{equation}
I_j = \{i : b_\rho(\Phi)_{ij} \ne 0 \}.
\end{equation}
\vskip 0.3cm

\begin{lem} 
Assume that $S(\Phi(D_\rho)) = S(D_\rho)$. 
Then, for each $j$,  
$$\lambda_i = \lambda_k \quad \text{for all } \quad i, k \in I_j.$$
\end{lem}
\vskip 0.3cm

\noindent{\it Proof}.  
Since $\eta$ is concave,  by Lemma 3.1 (2) we have  that 
\begin{eqnarray*}
\lefteqn{\sum_i \eta(\lambda_i) =  \sum_i \sum_j  b_\rho(\Phi)_{ij} \eta(\lambda_i) =  \sum_j \sum_i b_\rho(\Phi)_{ij} \eta(\lambda_i)} \\
 & \leqq&  \sum_j \eta(\sum_i \lambda_i b_\rho(\Phi)_{ij} ) = \sum_j \eta(\mu_j) = S(\Phi(D_\rho)) = S(\rho) = \sum_i \eta(\lambda_i)
\end{eqnarray*}
so that 
$\eta(\sum_i b_\rho(\Phi)_{ij} \lambda_i  ) = \sum_i b_\rho(\Phi)_{ij} \eta(\lambda_i)$.  

This  implies that  
$\lambda_i = \lambda_k$ for all $i, k \in I_j$ because $\eta$ is strictly concave.  
\qed
\vskip 0.3cm

Under the assumption that $S(\Phi(D_\rho)) = S(D_\rho)$, for $i \in I_j$ we denote the constant  $\lambda_i$ 
in  Lemma 3.2 by $\lambda^{(j)}$. 
Remark that each  $I_j$ is a non empty set because $b_\rho(\Phi)$ is a bistochastic matrix, and  
\begin{equation}
\lambda^{(j)} = \frac{\sum_{i \in I_j} \lambda_i}{|I_j|} = \lambda_k \quad \text{for all} \quad k \in I_j. 
\end{equation}
\vskip 0.3cm

\begin{thm}
Let $\rho$ be a state of $M_n(\mathbb {C})$  and 
let $\Phi : M_n(\mathbb {C}) \rightarrow  M_n(\mathbb {C})$ be a unital positive {\rm Tr}-preserving map. 
Then the following are equivalent{\rm :}
\smallskip

{\rm(i)} $S(\rho \circ \Phi^*) = S(\rho)$, i.e. $S(\Phi(D_\rho )) = S(D_\rho)$.  
\smallskip

{\rm(ii)} $\lambda  = \mu b_\rho(\Phi)^T$, 
where $x^T$ denotes the transpose of the matrix $x$. 
\vskip 0.3cm

{\rm(iii)} $\lambda_i = \mu_i$ for all $i = 1, \cdots, n$. 
\vskip 0.3cm

{\rm(iv)} $\Phi(D_\rho) =  u D_\rho u*$ for some unitary $u \in M_n(\mathbb {C})$. 
\vskip 0.3cm

{\rm (v)} $\Phi^* \Phi(D_\rho) = D_\rho$. 

\end{thm}
\vskip 0.3cm

\noindent{\it Proof.} 
(i) $\Rightarrow$ (ii): 
By Lemma 3.2,  we have    that  for each $j$, 
$$\mu_j = \sum_i \lambda_i  b_\rho(\Phi)_{ij} = \sum_{i \in I_j} \lambda_i  b_\rho(\Phi)_{ij} =  
 \lambda^{(j)} \sum_{i \in I_j}  b_\rho(\Phi)_{ij} =  \lambda^{(j)}.$$ 
If $b_\rho(\Phi)_{kj} \ne 0$, then $k \in I_j$ and 
$$
\sum_{j=1}^n \mu_j b_\rho(\Phi)_{kj}  = \sum_{j=1}^n \lambda^{(j)}  b_\rho(\Phi)_{kj} 
 = \sum_{j \in \{l : b_\rho(\Phi)_{kl} \ne 0 \}} \lambda^{(j)}  b_\rho(\Phi)_{kj} 
 =\lambda_k,
$$
this means that $\mu b_\rho(\Phi)^T = \lambda$.  
\smallskip

(ii) $\Rightarrow$ (iii): 
Remember the following fact in the majorization theory (for example, see \cite{Pe}): 
if $\lambda$ and $\mu$ are probability vectors such that $\lambda b = \mu$ for a bistochastic matrix $b$, then 
$\mu$ is majorized by $\lambda$, that is, 
$$\sum_{j=1}^k \mu_j  \leq \sum_{i=1}^k  \lambda_i \quad \text{for all} \quad  k = 1, \cdots, n.$$
The relation (ii) implies that $\lambda$ is majorized by $\mu$ and also Lemma 3.1 implies that 
$\mu$ is majorized by $\lambda$. Hence $\lambda_i = \mu_i$ for all $i = 1, \cdots, n$. 
\smallskip

(iii) $\Rightarrow  $ (iv):  Assume that $\lambda_i = \mu_i$ for all $i$. Then 
$$\Phi(D_\rho) = \sum_j \mu_j p_j = \sum_j \lambda_j u_{(\rho, \Phi)} e_j u_{(\rho, \Phi)}^* = u_{(\rho, \Phi)} \ D_\rho \ u_{(\rho, \Phi)}^* .$$ 

(iv) $\Rightarrow$  (v):  
It is sufficient to show that 
$ ||\Phi^* \Phi(D_\rho) - D_\rho ||_2 = 0$. 

Since $\Phi$ is unital, i.e.,  $\Phi(1) = 1$, it implies that  Tr $\Phi^* = $Tr. 
Also since the positive map $\Phi^*$ satisfies 
the Kadison-Schwartz inequality for the positive operator $\Phi(D_\rho)$,  
we have that 
\begin{eqnarray*}
\lefteqn{|| \Phi^* (\Phi(D_\rho)) ||_2^2 
 =   \text{Tr} (\ (\Phi^* (\Phi(D_\rho) ))^* (\Phi^* (\Phi(D_\rho))) \ ) }\\
 &\leq & \text{Tr} ( \Phi^* (\Phi(D_\rho) ^*  \Phi(D_\rho)) \ ) = \text{Tr} ((\Phi(D_\rho)) ^*  \Phi(D_\rho) ).  
\end{eqnarray*}
Hence  by the condition (iv) 
$$|| \Phi^* (\Phi(D_\rho) ) ||_2 \leq  \text{Tr} ((\Phi(D_\rho)) ^*  \Phi(D_\rho) ) 
= \text{Tr} ( u D_\rho ^* u^* u D_\rho u^*) ) = || D_\rho||_2$$ 
and 
$$<\Phi^* \Phi(D_\rho) , D_\rho>_{Tr} = \text{Tr} ((\Phi(D_\rho)) ^*  \Phi(D_\rho) ) = || D_\rho||_2^2.$$
These imply   that
\begin{eqnarray*}
\lefteqn{0 \leq ||\Phi^* \Phi(D_\rho) - D_\rho ||_2^2} \\
 &=&   \text{Tr} ( (\Phi^* \Phi(D_\rho) - D_\rho)^* (\Phi^* \Phi(D_\rho) - D_\rho) )\\
 &=& || \Phi^* \Phi(D_\rho) ||_2^2 - 2 <\Phi^* \Phi(D_\rho) , D_\rho> + ||D_\rho||_2^2 \\
 &\leq&  ||D_\rho||_2^2 -   ||D_\rho||_2^2 = 0
\end{eqnarray*}
so that  $\Phi^* \Phi(D_\rho) = D_\rho$. 
\smallskip

(v) $\Rightarrow$ (i):  Since $\Phi$ and $\Phi^*$ are unital positive {\rm Tr}-preserving, we have that 
$S(D_\rho) = S(\Phi^* \Phi (D_\rho)) \geq  S(\Phi (D_\rho)) \geq  S(D_\rho) $ 
so that $S(D_\rho) =  S(\Phi (D_\rho))$. 
\qed
\vskip 0.3cm 

\begin{rem}
Under the assumption that $\Phi$ is 2-positive, 
the corresponding  relation to (i) $\Leftrightarrow $ (v)  in  Theorem 3.3 is 
obtained for the discussion on  the relative entropy in \cite [Theorem 7.1]{HMPB} 
(cf. \cite{ZW} as an application of \cite{HMPB}).   

In our case, $\Phi$  is not necessary  to be 2-positive.  
\end{rem}

\vskip 0.3cm

\begin{ex}  Now we pick up the transpose mapping $\Phi: x \to x^T$ on $M_n(\mathbb{C})$. 
It is a typical example of unital Tr-preserving positive but not 2 positive  map. 
The $\Phi$ satisfies the  conditions in Theorem 3.3 for all state $\rho$. 

In fact, the $\Phi$ is a symmetry  as follows: 
\begin{eqnarray*}
\lefteqn{ <\Phi^*(x), y> = <x, \Phi(y)> = {\rm Tr}(\Phi(y)^*x) = {\rm Tr}((y^T)^*x) } \\
 &=& \sum_{i,j = 1}^n \overline {y_{i,j}} x_{j,i} = \sum_{i,j = 1}^n \overline{y_{j,i}} x_{i,j} = {\rm Tr}(y^* x^T) 
      = {\rm Tr}(y^* \Phi(x)) \\
 &=&  <\Phi(x), y> \quad \text{for all} \quad x = (x_{ij}), \ y = (y_{ij}).  
\end{eqnarray*}
Hence $\Phi^* \Phi$ is the identity map on $M_n(\mathbb{C})$ so that (v) in Theorem 3.3 is trivial. 
\end{ex}

\vskip 0.3cm

\begin{rem} 
If the state $\rho$ in Theorem 3.3 is the normalized trace $\text{Tr} / n$, then  $D_\rho = \text{I}_n / n$  
so that the statements (i) - (v) are all trivial for every $\Phi$. 
\end{rem}

\begin{rem} 
The statements (i) and (iv)  in Theorem 3.3 are  extended versions of  the statements 1 and 4 
in \cite [Theorem 4.5]{PP} respectively. 
The statement (ii)  in Theorem 3.3 is  nothing else but 
$ \lambda  = \lambda b_\rho(\Phi) b_\rho(\Phi)^T,$  
which corresponds with the statement 3 in \cite [Theorem 4.5]{PP}. 
\end{rem}
\vskip 0.3cm

\begin{cor}
{\rm (1)} If  $S(\Phi(D_\rho)) = S(\rho)$, then 
$$<\Phi(D_\rho), \Phi(e_k) > =  <D_\rho, e_k> \ \text{for all} \ k.$$

{\rm (2)} If $E_B$ is the onditional expectation of $M_n(\mathbb {C})$ onto $B$, 
then 
$$S(E_B(D_\rho))  = S(D_\rho) \quad \text{if and only if} \quad D_\rho \in B.$$

\end{cor}
\smallskip

\noindent{\it Proof.} 
(1): The assumption and (ii) of Theorem 3.3 imply 
that for each $k$ 
\begin{eqnarray*}
\lefteqn{\text{Tr} ( \Phi^* (\Phi(D_\rho)) e_k) =  \text{Tr} (\Phi(D_\rho) \Phi(e_k)) = \text{Tr} ( \sum_{j=1}^n \mu_j p_j\Phi(e_k) ) } \\
 &=&\sum_{j=1}^n \mu_j \text{Tr} ( p_j\Phi(e_k) ) = \sum_{j=1}^n \mu_j b_\rho(\Phi)_{kj} = \lambda_k, 
\end{eqnarray*}
so that 
$<\Phi(D_\rho), \Phi(e_k) > = \text{Tr} ( \Phi^* \Phi(D_\rho) e_k) = \lambda_k =  <\rho, e_k>$ for all $k$.
\smallskip

(2): A conditional expectation $E$ satisfies that $E^*E = E$. 
Hence  by (i) $\Leftrightarrow $ (v)  in Theorem 3.3,  we have that 
$S(E_B(D_\rho))  = S(D_\rho)$ if and only if $D_\rho = E_B^*E_B(D_\rho) = E_B(D_\rho)$  which means 
that $D_\rho \in B$. 
\qed
\vskip 0.3cm

\subsection{\bf Relations among various entropies} 
In this section, we discuss about various kinds of entropy under the same notations with in Section 3.1: 
for  a state $\rho$ and a positive unital Tr-preserving map $\Phi$ of the algebra $M_n(\mathbb {C})$, 
the probability vector 
$\lambda = (\lambda_1, \cdots, \lambda_n )$ 
(resp. $\mu = ( \mu_1, \cdots, \mu_n )$ ) is given by eigenvalues of $D_\rho$ (resp. $\Phi(D_\rho)$) 
whose coresponding minimal projections are 
$\{ e_1, \cdots, e_n \}$ (resp. $\{ p_1, \cdots, p_n \}$ ), 
and $A$ (resp. $B$) is the  subalgebra generated by $\{e_i\}_{i=1}^n$ (resp. $\{p_j\}_{j=1}^n$ ). 
\vskip 0.3cm

{\bf 3.2.1.}  
We apply the notion of the weighted entropy for a bistochastic matrix defined in \cite{ZSKS} 
to our bistochastic matrix $b_\rho(\Phi)$.   
We let 
$$H^\lambda(b_\rho(\Phi)) = \sum_{i=1}^n  \lambda_i \sum_{j=1}^n \eta(b_\rho(\Phi)_{ij}) \quad \text{and } \quad 
H_\mu(b_\rho(\Phi)) = \sum_{j=1}^n \mu_j \sum_{i=1}^n \eta(b_\rho(\Phi)_{ij}).$$ 
\smallskip

This is well defined, i.e. 
the values $H^\lambda(b_\rho(\Phi))$ and $H_\mu(b_\rho(\Phi))$ depend on only the pair $\{\rho, \Phi\}$. 
In fact, assume that 
$D_\rho = \sum_i \lambda_i f_i$ and  $\Phi(D_\rho) =  \sum_j \mu_j q_i$ 
for minimal projections $\{f_i\}_{i=1}^n$ and $\{q_j\}_{j=1}^n$ which  are not always 
same as  $\{e_i\}_{i=1}^n$ and $\{p_j\}_{j=1}^n$.  
Then there are permutations  $\sigma$ and $\pi$ of $\{1, \cdots, n\}$ such that 
$f_i = e_{\sigma(i)}$ and $q_i = p_{\pi(i)}$ for all $i$. 
Remark that $\lambda_i = \lambda_{\sigma(i)}$ for all $i$, then   
\begin{align*}
\lefteqn{\sum_{i=1}^n \lambda_i \sum_{j=1}^n \eta( \text{Tr}(\Phi(f_i)q_j) ) 
           = \sum_{i=1}^n \lambda_i \sum_{j=1}^n \eta( \text{Tr}(\Phi(e_{\sigma(i)}) p_{\pi(j)}) ) } \\
 &= \sum_{i=1}^n \lambda_{\sigma(i)} \sum_{j=1}^n \eta( \text{Tr}(\Phi(e_{\sigma(i)}) p_{\pi(j)}) )  
 = \sum_{i=1}^n \lambda_i \sum_{j=1}^n \eta( \text{Tr}(\Phi(e_i) p_j) ).  
\end{align*}
Hence the value $H^\lambda(b_\rho(\Phi))$ does not depend on the choice of  minimal projections. 
Similarly it holds for $H_\mu(b_\rho(\Phi))$ by that $\mu_j = \mu_{\pi(j)}$ for all $j$. 
\vskip 0.3cm

{\bf 3.2.2.}  
Since $\Phi$ is positive unital Tr-preserving map,  $\Phi(e_i)$ and $\Phi^*(p_j)$ are  density matrices for all $i,j$. 
We put  $S_\rho(\Phi)$ and  $S^\rho(\Phi^*)$ as the followings: 
$$ S_\rho(\Phi) = \sum_{i=1}^n \lambda_i S(\Phi(e_i)) \quad \text{and} \quad 
S^\rho(\Phi^*) = \sum_{j=1}^n \mu_j S(\Phi^*(p_j)). $$
\vskip 0.3cm

Similarly to the case of $H^\lambda(b_\rho(\Phi))$ and $H_\mu(b_\rho(\Phi))$, 
the values $S_\rho(\Phi)$ and $S^\rho(\Phi^*)$ are  uniquely determined by the pair $\{\rho, \Phi\}$.  
\vskip 0.3cm

\begin{prop}
For the conditional expectation $E_B$ (resp. $E_A$) onto  $B$ (resp. $A$), 
the followings hold: 
\begin{enumerate}
 \item $E_B(\Phi(e_i)) = \sum_{j=1}^n b_\rho(\Phi)_{ij} p_j$,   $ E_A(\Phi^*(p_j)) = \sum_{i=1}^n b_\rho(\Phi)_{ij} e_i$, 
 \item $H^\lambda(b_\rho(\Phi)) = \sum_i \lambda_i S(E_B(\Phi(e_i)))$, 
  $H_\mu(b_\rho(\Phi)) = \sum_j \mu_j S(E_A(\Phi^*(p_j))).$ 
\end{enumerate}
\end{prop}
 \smallskip

\noindent{\it Proof.} 
1: It is clear that  $E_B(\Phi(e_i))$  is a density matrix for all $i$. 
The $E_B$ is given by 
$E_B(x) = \sum_{j=1}^n  {\rm Tr}(p_jx)  p_j$ for all $x \in M_n(\mathbb {C})$. 
Hence we have  the form  for $E_B(\Phi(e_i))$ and similarly  for $ E_A(\Phi^*(p_j))$. 
\smallskip 

2:  By the definition, we have that   
$H^\lambda(b_\rho(\Phi)) 
= \sum_{i=1}^n \lambda_i \sum_{j=1}^n \eta(b_\rho(\Phi)_{ij})
= \sum_{i=1}^n \lambda_i S(E_B(\Phi(e_i)))$ 
and similarly  the statements about $H_\mu(b_\rho(\Phi))$. 
\qed
\vskip 0.3cm

For a probability vector $\lambda = (\lambda_1, \cdots, \lambda_n )$, we set 
\begin{equation}
J_\lambda = \{k ; \lambda_k \ne 0 \}. 
\end{equation}

\begin{thm} 
Let $\rho$ be a state  of $M_n(\mathbb {C})$,  and let $\Phi$ be a unital positive {\rm Tr}-preserving map on $M_n(\mathbb {C})$. 
Then the following statements hold: 
\smallskip 

{\rm (1)} 
\begin{align*}
S_\rho(\Phi) \leq H^\lambda(b_\rho(\Phi)) &\le S(\rho \circ \Phi^*)    \le  S(\rho) + S_\rho(\Phi) \\
 &\hskip 1.1cm \vee \Vert \\
 &\hskip 1.1cm   S(\rho)
\end{align*}

{\rm (2)}  $S_\rho(\Phi) =  H^\lambda(b_\rho(\Phi))$ if and only if $\Phi(e_i) \in B$ for all $i \in   J_\lambda$:  
$$\Phi(e_i) = \sum_{j=1}^n b_\rho(\Phi)_{ij} p_j, \quad \text{for all} \ i \in  J_\lambda.$$ 
\smallskip

{\rm (3)} 
$H^\lambda(b_\rho(\Phi)) =  S(\rho \circ \Phi^*)$ if and only if $b_\rho(\Phi)_{ij} = \mu_j$ for all  $i \in  J_\lambda$ and $j${\rm :}  
$$ \Phi(D_\rho) = \sum_{j=1}^n b_\rho(\Phi)_{ij} p_j.$$ 
\smallskip

{\rm (4) } 
$ S_\rho(\Phi) = S(\rho \circ \Phi^*)$ if and only if  $\Phi(D_\rho) = \Phi(e_i)$ for every  $i \in J_\lambda$.   
\smallskip

{\rm (5) } 
$S(\rho \circ \Phi^*) = S(\rho) + S_\rho(\Phi) $ if and only if  the $\rho$ is a pure state. 
\end{thm} 
\smallskip

\noindent{\it Proof.} 
(1):  
Since $E_B$ is a unital positive Tr-preserving mapping, we have that 
$S(\Phi(e_i)) \leq S(E_B(\Phi(e_i)))$ for all $i$,  
and by Proposition 3.9 
$$S_\rho(\Phi) = \sum_{i=1}^n \lambda_i S(\Phi(e_i)) \leq \sum_{i=1}^n \lambda_i S(E_B(\Phi(e_i))) 
= H^\lambda(b_\rho(\Phi)).$$ 
On the other hand, $\eta$ is concave and $\lambda b_\rho(\Phi) = \mu$, we have that 
$$H^\lambda(b_\rho(\Phi)) = \sum_j \sum_i \lambda_i \eta( b_\rho(\Phi)_{ij})
\leqq \sum_j \eta( \sum_i \lambda_i b_\rho(\Phi)_{ij})  = S(\Phi(D_\rho)).$$
The inequality $S(\Phi(D_\rho)) \le S(D_\rho) + S_\rho(\Phi)$ is induced by the inequality for 
convex conbinations of density operators (cf. \cite{NS, OP, Pe}) but 
we denote a computation  which we need  below to prove (5).  

Since $\Phi$ is Tr-preserving,  we have that 
\begin{eqnarray*} 
\lefteqn{S(\Phi(D_\rho)) = S(\sum_i \lambda_i \Phi(e_i) ) 
 \leq  \sum_i  S(\lambda_i \Phi(e_i) ) = \sum_i  \text{Tr} \left(  \eta(\lambda_i \Phi(e_i) ) \right) }\\
 &=& \sum_i \text{Tr} ( \eta( \lambda_i) \Phi(e_i)  + \lambda_i \eta( \Phi(e_i) ) \ ) 
 =  S(D_\rho) + S_\rho(\Phi). 
\end{eqnarray*}
The relation that  $S(\rho) \le S(\rho \circ \Phi^*)$ is known, and  
the all relations in (1) hold.
\vskip 0.3cm

(2):  If $S_\rho(\Phi) = H^\lambda(b_\rho(\Phi))$, then 
$\sum_i \lambda_i S(\Phi(e_i)) = \sum_i \lambda_i S( E_B(\Phi(e_i)))$. 
Since $S(\Phi(e_i)) \le S( E_B(\Phi(e_i)))$ for all $i$, this impies  that  
$\lambda_i  S(\Phi(e_i)) = \lambda_i S( E_B(\Phi(e_i)))$ for all $i$ so that 
$S(\Phi(e_i)) = S(E_B(\Phi(e_i)))$ for all $i \in J_\lambda$. 

Now we apply Theorem 3.3 to the pair $\{\Phi(e_i), E_B\}$ for each $i \in J_\lambda$. 
Then by the property of the conditional expectation that $E_B = E_B^* E_B$ 
$$E_B(\Phi(e_i)) = E_B^* E_B(\Phi(e_i)) = \Phi(e_i) \quad \text{for each} \quad  i \in J_\lambda.$$
This means that  $\Phi(e_i) \in B$ for all $i \in J_\lambda$.  

Conversely assume that $\Phi(e_i) \in B$ for all $i \in J_\lambda$. 
Then 
$$S_\rho(\Phi) =  \sum_{i \in J_\lambda} \lambda_i S(\Phi(e_i)) 
= \sum_{i \in J_\lambda} \lambda_i S(E_B(\Phi(e_i))) = H^\lambda(b_\rho(\Phi)).$$

We remark that $\Phi(e_i) \in B$ if and only if $\Phi(e_i) = \sum_{j=1}^n b_\rho(\Phi)_{ij} p_j$. 
In fact, if $\Phi(e_i) \in B$, then  there exists $\{ \alpha_{ij} \}_j \subset \mathbb{C}$ such that 
$\Phi(e_i) = \sum_j \alpha_{ij} p_j$, which implies that 
$\alpha_{ij} = \rm{Tr}(\Phi(e_i) p_j) = b_\rho(\Phi)_{ij}$  so that 
$\Phi(e_i) = \sum_{j=1}^n b_\rho(\Phi)_{ij} p_j$.  
It is obvious that $\Phi(e_i) = \sum_{j=1}^n b_\rho(\Phi)_{ij} p_j \in B$. 
\smallskip 

(3):  
Assume that $H^\lambda(b_\rho(\Phi)) = S(\rho \circ \Phi^*)$ (so that $= S(\Phi(D_\rho)) $). 
Then 
\begin{eqnarray*}
\lefteqn{\sum_j \sum_i \lambda_i \eta(b_\rho(\Phi)_{ij})  = H^\lambda(b_\rho(\Phi)) = S(\Phi(D_\rho))} \\
 &=& \sum_j \eta (\mu_j) = \sum_j \eta (\sum_i \lambda_i b_\rho(\Phi)_{ij} ).
\end{eqnarray*}
This implies  that 
$$\sum_i \lambda_i \eta(b_\rho(\Phi)_{ij})  = \eta (\sum_i \lambda_i b_\rho(\Phi)_{ij} ) \quad \text{for all } \ j$$
because in general $\sum_i \lambda_i \eta(b_\rho(\Phi)_{ij})  \le \eta (\sum_i \lambda_i b_\rho(\Phi)_{ij} )$ 
for all $j$ by the fact that $\eta$ is concave. 
Moreover  $\eta$ is strictly concave, so that  for each $j$    
$$b_\rho(\Phi)_{ij} = b_\rho(\Phi)_{i'j} \quad \text{for all} \quad i, i' \in   J_\lambda.$$
The $\lambda = (\lambda_1, \cdots, \lambda_n)$ is a probability vector 
so that  there exists a $k$ with $\lambda_k \ne 0$, which we fix. 
If $\lambda_i \ne 0$, then $b_\rho(\Phi)_{ij} = b_\rho(\Phi)_{kj}$ for all $j$. 
Hence 
$$\mu_j = \sum_i \lambda_i b_\rho(\Phi)_{ij} = \sum_{i \in J_\lambda} \lambda_i b_\rho(\Phi)_{ij} 
= \sum_{i \in J_\lambda} \lambda_i b_\rho(\Phi)_{kj} = b_\rho(\Phi)_{kj}, \ \text{for all} \ j.$$ 

Conversely assume that for each $j$,  $b_\rho(\Phi)_{ij}  = \mu_j$ for all $i \in J_\lambda$.  
Then 
$$H^\lambda(b_\rho(\Phi)) = \sum_{i \in J_\lambda} \lambda_i \sum_j  \eta(b_\rho(\Phi)_{ij}) = 
\sum_{i \in J_\lambda} \lambda_i \sum_j  \eta(\mu_j) =  S(\Phi(D_\rho)).$$

(4):  
If $S_\rho(\Phi) = S(\rho  \circ  \Phi^*)$, then 
$\text{Tr}(\sum_i \lambda_i \eta(\Phi(e_i))) = \text{Tr} (\eta (\sum_i \lambda_i \Phi(e_i))$. 
This implies that 
$\sum_{i \in J_\lambda} \lambda_i \eta(\Phi(e_i)) = \eta (\sum_{i \in J_\lambda} \lambda_i \Phi(e_i))$ 
because $\eta$ is operator-concave and $\text{Tr}$ is faithful. 
Moreover $\eta$ is strictly operator-concave. 
Hence $\Phi(e_i) = \Phi(e_j)$ for all $i, j \in J_\lambda$ and we have   that 
$$\Phi(D_\rho) = \sum_{i \in J_\lambda} \lambda_i \Phi(e_i) = \Phi(e_i) \ \text{for all} \ i \in J_\lambda.$$ 

Conversely assume that $\Phi(D_\rho) = \Phi(e_i)$ for every $i \in J_\lambda$. 
Then $\Phi(e_k) = \Phi(e_i)$ for all $k, i \in J_\lambda$. 
Let us fix an $i$ in $J_\lambda$. 
Then for all $j = 1, \cdots, n$ 
\begin{eqnarray*}
\lefteqn{\mu_j = \sum_{k} \lambda_k b_\rho(\Phi)_{kj} =  \sum_{k \in J_\lambda}  \lambda_k b_\rho(\Phi)_{kj} 
 = \sum_{k \in J_\lambda} \lambda_k {\rm Tr}(\Phi(e_k)p_j)} \\
 &=&  \sum_{k \in J_\lambda} \lambda_k {\rm Tr}(\Phi(e_i)p_j) = {\rm Tr}(\Phi(e_i)p_j) = b_\rho(\Phi)_{ij}.
\end{eqnarray*}
Hence by (3) we have  that  $H^\lambda(b_\rho(\Phi)) = S(\rho \circ \Phi^*)$ . 

On the other hand, the assumption that $\Phi(D_\rho) = \Phi(e_i)$ for every $i \in J_\lambda$ implies clearly 
that $\Phi(e_i) \in B$ for all $i \in J_\lambda$ so that 
$S_\rho(\Phi) = H^\lambda(b_\rho(\Phi))$ by (2). 
Thus we have that  
$S_\rho(\Phi) = S(\rho \circ \Phi^*)$. 
\smallskip

(5) 
Assume that $S(\rho \circ \Phi^*) = S(\rho) + S_\rho(\Phi)$. 
By the computation in (1) for $S(\rho \circ \Phi^*) \le S(\rho) + S_\rho(\Phi)$,  
we have that $\text{Tr} (\eta(\sum_i \lambda_i \Phi(e_i) ) )  = \text{Tr} (\sum_i \eta(\lambda_i \Phi(e_i) ))$. 
On the other hand, 
$\eta(\sum_i \lambda_i \Phi(e_i) ) \le \sum_i \eta(\lambda_i \Phi(e_i) )$ 
because $\eta$ is operator concave.  
These two relations imply that 
$\eta(\sum_i \lambda_i \Phi(e_i) )   = \sum_i \eta(\lambda_i \Phi(e_i) )$ 
by the faithfulness of  $\text{Tr}$.   
Again by the fact that $\eta$ is strictly operator concave,   we have  that 
$\Phi(e_i) = \Phi(e_j)$ for all $i,j \in J_\lambda$ so that $\Phi(D_\rho) = \Phi(e_i)$ for all $i \in J_\lambda$. 
By combinig with (4),  this  implies that  $S_\rho(\Phi) = S(\rho \circ \Phi^*)$ 
so that by the assumption $S(\rho) = 0$  and so  $\rho$ is a pure state,  

Conversely if $\rho$ is a pure state, then $S(\rho) = 0$. 
Hence by (1) $S_\rho(\Phi) \le S(\rho \circ \Phi^*) \le S_\rho(\Phi)$   
so that  $S(\rho \circ \Phi^*) = S_\rho(\Phi) = S(\rho) + S_\rho(\Phi)$. 
\qed
\smallskip 

\begin{cor} 
Assume that all eigenvalues of $D_\rho$ are  nonzero and that 
$H^\lambda(b_\rho(\Phi)) = S(\rho \circ \Phi^*)$. 
Then for all $\Phi$ the state $\rho \circ \Phi^*$ is the canonical tracial state, i.e. ${\rm Tr} / n$ 
and  so that $H^\lambda(b_\rho(\Phi)) = S(\rho \circ \Phi^*) = \log n$. 
\end{cor}
\smallskip 

\noindent{\it Proof.} 
If  $\lambda_i \ne 0$ for all $i$ and if $H^\lambda(b_\rho(\Phi)) = S(\Phi(D_\rho)) $, 
then by (3) of Theorem 3.10  $\mu_j = b_\rho(\Phi)_{ij}$ for all $i, j$ so that $n\mu_j = \sum_i b_\rho(\Phi)_{ij} = 1$  
for all $j$. Hence $\mu_j = 1/n$ for all $j$.    
This means that $\Phi(D_\rho)$ is the density matrix of the normalized trace ${\rm Tr} / n$ so that $S(\Phi(D_\rho)) = \log n$. 
\qed
\smallskip 

\begin{rem} 
A bistochastic matrix $b$ is said to be {\it unistochastic } if 
it is induced from some unitary matrix $u$ by the method that $b_{i, j} = |u_{i, j}|^2$ for all $i, j = 1, \cdots, n$. 
An $n \times n$ unitary matrix $u$ is called a {\it Hadamard matrix} if 
$|u_{i, j}| = 1 / \sqrt n$ for all $i, j = 1, \cdots, n$. 

The  above corollary shows that if all eigenvalues of $D_\rho$ are non-zero 
and if $H^\lambda(b_\rho(\Phi)) = S(\Phi(D_\rho))$ then $b_\rho(\Phi)$ is a 
unistochastic matrix induced from a Hadamard matrix. 
\end{rem}
\smallskip 

\begin{ex} 
(1) Assume that $\rho$ is a pure state. Then  it is clear that 
$$S(\rho) = 0 \quad \text{and} \quad S_\rho(\Phi) = H^\lambda(b_\rho(\Phi)) = S(\rho \circ \Phi^*)$$
for every positive unital Tr-preserving map $\Phi$. 
\smallskip 

Furthermore, for any given value $s$ with $0 \le s \le \log n$, there exists a positive unital Tr-preserving map $\Phi$ such that 
$S_\rho(\Phi) = s$. 

In fact, we may assume that $D_\rho = e_1$. Since  $\{e_1,  \cdots, e_n\}$ are minimal mutually orthogonal projections, 
there exists a family of  partial isometries $\{ v_{j i} : j, i  = 1, \cdots, n\}$ such that 
$v_{j i} v_{j i}^* = e_j$ and $v_{j i}^* v_{j i}= e_i$. 
Choose an $n$-tupple of numbers $\mu_j$ such that 
$\mu_j \ge 0$,  $\sum_{j=1}^n \mu_j = 1$ and $\sum_{j=1}^n \eta(\mu_j) = s$. 
We define the map $\Phi$  by 
$$\Phi(x) = \sum_{k = 1}^n \sum_{j=1}^n \mu_{j + k-1} v_{jk} x v_{jk}^* \ ( \bmod \ n), \quad x \in M_n(\mathbb {C}).$$
Then $\Phi$ satisfies the conditions and   
$$S_\rho(\Phi)  = S(\Phi(e_1)) = {\rm Tr} (\eta (\sum_j \mu_j e_j) ) = \sum_j \eta(\mu_j) = s.$$
\smallskip 

(2) 
If $\Phi$ is a *-isomorphism, then for each state  $\rho$,  
$$  S_\rho(\Phi) = H^\lambda(b_\rho(\Phi)) = 0 \quad  \text{and} \quad  S(\rho \circ \Phi^*) =  S(\rho).$$
 
In fact, $\Phi(e_i)$ is a minimal projection for all $i$ and 
the  set $\{\lambda_i\}_i$  coincides with the  set $\{\mu_i\}_i$. 
Hence 
$S_\rho(\Phi) = \sum_{i=1}^n \lambda_i S(\Phi(e_i)) = 0$ and 
${\rm Tr}(\Phi(e_i) p_j))$ is either $0$ or $1$ for all $i, j$
so that $H^\lambda(b_\rho(\Phi)) = 0$. 
Obviously  $S(\Phi(D_\rho)) = S(D_\rho)$. 
\smallskip 

(3) If $\Phi$ is  a map $M_n(\mathbb {C}) \to \mathbb{C}1_M$,  then  for every state  $\rho$ 
$$S_\rho(\Phi) = H^\lambda(b_\rho(\Phi)) =  S(\rho \circ \Phi^*)  = \log n.$$

In fact, each $\Phi(e_i)$ is $\alpha_i 1_M$ for some $\alpha_i \in \mathbb{C}$. 
Hence $1 = \rm{Tr}(e_i) = \rm{Tr}(\Phi(e_i) ) = \alpha_i \rm{Tr}(\Phi(1_M)) =  \alpha_i  n$, i.e., 
$\Phi(e_i) = \frac 1n 1_M$. 
This implies that 
$S_\rho(\Phi) = \sum_i \lambda_i S(1_M/n ) = \text{Tr} (\eta(1_M / n) ) = \log n$. 
In general  
$S_\rho(\Phi) \le  H^\lambda(b_\rho(\Phi)) \le  S(\rho \circ \Phi^*) \le \log n$ for every state $\rho$.  
Hence we have the conclusion.   
\smallskip

(4) A typical counter example of $\Phi$ for Theorem 3.10 (4)  is the transpose mapping $\Phi(x) = x^T$, where $\rho$ is not a pure state.  
\smallskip

More general examples are given as follows: 
Let $D = \sum_i \lambda_i e_i$ be a given density matrix. Let $\{p_j\}_{j=1}^n$ be mutually orthogonal minimal  projections. 
Then we have a family of partially isometries $\{v_{ij} \}_{ij} $ such that $v_{ij}^* v_{ij} = e_j$ and $v_{ij} v_{ij}^* = p_j$. 
Let $a = [a_{ij}]$ be a bistochastic matrix and let 
$$\Phi(x) = \sum_{i,j} a_{ij} v_{ij} x v_{ij}^*, \quad (x \in M_n(\mathbb {C})).$$
Then $\Phi$ is a unital positive Tr-preserving map and 
$$\Phi(D) = \sum_i (\sum_j a_{ij} \lambda_j) p_i \quad \text{and}  \quad \Phi(e_i) = \sum_j a_{ji} p_j \quad \text{for all} \ i.$$

Hence $\Phi(e_i) \in B$ for all $i$, that is the condition in Theorem 3.10 (2). 
\smallskip

Also we can choose  bistochastic matrices $a = [a_{ij}]$, one of which induces $\Phi$ satisfying the condition  (4) in Theorem 3.10 and 
the other of which induces $\Phi$ not satisfying the condition  (4) in Theorem 3.10. 
\end{ex}
\vskip 0.3cm

At the last, we give the following characterization: 
\vskip 0.3cm

\begin{thm} 
Let $\rho$ be a state  of $M_n(\mathbb {C})$,  and let $\Phi$ be a unital positive {\rm Tr}-preserving map on $M_n(\mathbb {C})$. 
Then the following  conditions  are equivalent: 

\rm 0) $H^\lambda(b_\rho(\Phi)) = 0$, 

\rm 1) for each  $i \in J_\lambda$, there exists a unique $j(i)$ such that 
$$\lambda_i = \mu_{j(i)} \quad \text{and} \quad \Phi(e_i) = p_{j(i)},$$

\rm 2) $S(\rho) = S(\rho \circ \Phi^*)$,

\rm 3) there exists a unitary $u $ such that $\Phi(D_\rho) = u D_\rho u^*$, 

\rm 4) $ \Phi^* \Phi(D_\rho) = D_\rho$.  
\end{thm}
\vskip 0.3cm

\noindent{\it Proof.} 0) $\Rightarrow $ 1):  
Assume that $H^\lambda(b_\rho(\Phi)) = 0$. Then 
$\sum_{i \in J_\lambda} \lambda_i \sum_j  \eta(b_\rho(\Phi)_{ij}) = 0$ 
which implies that, for each $i \in J_\lambda$, $b_\rho(\Phi)_{ij}$ is either $0$ or $1$ for all $j$. 

Let us fix an $i \in J_\lambda$.  Then we have a  unique $j(i)$ with 
$ b_\rho(\Phi)_{ij} =  \delta_{j, j(i)}$ by  that $\sum_j b_\rho(\Phi)_{ij} = 1$. 
Forthermore since $\sum_k b_\rho(\Phi)_{kj(i)} = 1$ it holds that 
$ b_\rho(\Phi)_{kj(i)} =  \delta_{k, i}$. 
As a consequence, for this $j(i)$, we have  that
$$b_\rho(\Phi)_{k j(i)} = \delta_{k, i} \quad \text{for all } \quad  k = 1, \cdots, n .$$

Remember  that 
$\lambda b_\rho(\Phi) = \mu$ by Lemma 3.1. 
This implies that 
$\mu_{j(i)} = \lambda_i$ for all $i \in J_\lambda$ 
so that 
$\Phi(e_i) = p_{j(i)}$ for all $i \in J_\lambda$ 
because $\sum_i  \lambda_i \Phi(e_i) = \Phi(D_\rho) = \sum_i p_{j(i)} p_{j(i)}$. 
\smallskip 

1) $\Rightarrow $ 0):  For all $i \in J_\lambda$,  
we have that  $E_B(\Phi(e_i)) = E_B(p_{j(i)}) = p_{j(i)}$ by the assumption,  so that 
$E_B(\Phi(e_i))$ is a minimal projection. 
This implies that $S(E_B(\Phi(e_i))) = 0$ for all $i \in J_\lambda$. 
Hence by  Proposition 3.9, we have that 
$$H^\lambda(b_\rho(\Phi)) = \sum_{i=1}^n \lambda_i S(E_B(\Phi(e_i)) = \sum_{i \in J_\lambda} \lambda_i S(E_B(\Phi(e_i)) = 0.$$

1) $\Rightarrow $ 2): Since $\Phi(D_\rho) = \sum_i \lambda_i \Phi(e_i) =  \sum_{i \in J_\lambda} \lambda_i \Phi(e_i)$, 
by combining the fact   $\Phi(e_i) = p_{j(i)}$ for all $i \in J_\lambda$, 
we have that $\{\lambda_i; i \in J_\lambda \}$ is the all non-zero eigenvalues of $\Phi(D_\rho)$ 
 which are of course the all eigenvalues of $D_\rho$. 
\smallskip 

It is obvious that 3) $\Rightarrow $ 1),  and  we have the conclusion by Theorem 3.3.
\qed

\end{document}